
\magnification=\magstep1
\pageno=0
\overfullrule=0pt
{1.5}

\font\eightrm=cmr8


\def\real{\mathop{{\rm I}\kern-.2em\hbox{\rm R}}\nolimits}
\def\footnoterule{\kern-3pt \hrule width 2truein \kern 2.6pt}
\def\bR{{\bf R}}  \def\bP{{\bf P}}   \def\bS{{\bf S}}
\def\bZ{{\bf Z}}  \def\bb{{\bf b}}
  \def\ba{{\bf a}}
\def\bPhi{{\bf \Phi}}  \def\bc{{\bf c}}  \def\bp{{\bf p}}  \def\ba{{\bf a}}
\def\bdelta{\mathop{\delta\kern-.45em\hbox{$\delta$}}\nolimits}
\def\SUBSECTION#1\par{\vskip0pt plus.2\vsize\penalty-75
        \vskip0pt plus-.2\vsize\bigskip\bigskip
        \leftline{\bf #1}\nobreak\smallskip\noindent}
\def\bbeta{\mathop{\beta\kern-.55em\hbox{$\beta$}}\nolimits}
\def\bp{{\bf p}}
\def\sqr#1#2{{\vcenter{\vbox{\hrule height.#2pt
               \hbox{\vrule width.#2pt height#1pt \kern#1pt
                \vrule width.#2pt}
              \hrule height.#2pt}}}}
 
\baselineskip=18pt

\centerline{Some Remarks on Multiresolution Analyses}
\centerline{Containing Compactly Supported Functions}

\bigskip
\centerline{V.\ Dobri{\' c}, R.\ F.\ Gundy, and P.\ Hitczenko}

\centerline{Lehigh University, Rutgers University, and North Carolina
State University}
\bigskip
\bigskip
\bigskip

\centerline{Rutgers University Technical Report \#97-003}
\midinsert
\vskip 2in
\endinsert



\eject

{\bf 1. Introduction.}
In [7], P.\ G.\ Lemari{\'e} proved that if a multiresolution analysis
contains a compactly supported function, then it contains a minimal
(pre)scale function. That is, a function $\phi(x)$ of compact support
such that

(1) the integer translates, $\phi(x - k),\,\, k\in\bZ$, are a Riesz basis for
the space $V_0$;

(2) every function in $V_0$ that is compactly supported may be written as a
finite linear combination of translates of $\phi$.

\noindent
The most basic examples of minimal scale functions of this type are the
$B$-splines, and the compactly supported scale functions constructed by
I.\ Daubechies [3].
An important property of these minimal-scale functions was first proved by
Y.\ Meyer [10] for Daubechies' functions, and subsequently 
stated by Lemari{\'e} [7] in the general case:
{\it The translates of $\phi$, restricted to the unit interval, form a
linearly independent set.}

The purpose of this note is to prove 
the following stronger version of the above for a
minimal scale function $\phi$ that is continuous on the unit interval: the
translates of $\phi$ are linearly independent over any subset of
positive measure contained in the unit interval. The stronger version is of
interest because it may be used to obtain a local convergence theorem for
multiresolution analyses with continuous minimal scale functions. The first
version of this type of local convergence theorem was proved by Gundy and
Kazarian [2] for a class of wavelet expansions that include the spline
wavelets. Their theorem assumed a regularity condition (condition (M-Z) of
[2]). It turns out that this regularity condition is, in fact, a
property of all multiresolution analyses with continuous minimal scale
functions, as a consequence of the above strong linear independence of these
functions.
\bigskip

{\bf 2. Notation.}
We suppose that a multiresolution analysis is given. That is, we have a scale
of subspaces of $L^2(\bR)$, 
$V_j,\,\, j\in \bZ$, such that $V_j \subset V_{j+1}$, and
$f(\cdot)\in V_j$ iff $f(2^{-j}\cdot)\in V_0$.
Furthermore, we are given a function $\phi\in V_0$ such that the integer
translates $\phi(\cdot - k)$ form a Riesz basis for $V_0$: any function
$f(\cdot)\in V_0$ has a representation
$$f(x) = \sum a_k \phi(x-k)$$
with
$$\sum a_w^2 \cong \Vert f\Vert_2^2.$$
If $\phi$ has a nonzero integral, then it follows that the increasing
sequence of subspaces exhausts $L^2(\bR)$.  (See [5, Chapter 2].) 
Let $P_j$ be the orthogonal projection operator from $L^2$ onto $V_j$.

Now we impose another restriction on the multiresolution analysis.  
We require that the space $V_0$ contain a nontrivial continuous function 
that is compactly supported.  (The continuity restriction omits the
Haar case, for which the theorems below are obviously true. However,
a more tedious formulation of assumptions, designed to include this
case, hardly seems worthwhile.  Our class of multiresolution analyses
does include the spline wavelets, the compactly supported Daubechies
wavelets, and those obtained from these classes by integration, as 
indicated in Lemari{\'e} [7].)
With this additional assumption, 
the techniques of Lemari{\'e} [7] may be used to show that there exists a
minimally supported, real-valued, continuous function $\phi\in V_0$ such
that every compactly supported function in $V_0$ admits a representation
as a finite linear combination of integer translates of $\phi$.
If we agree to normalize $\phi$ by setting its integral equal to one,
then $\phi$ is unique, up to integer translates.

\bigskip

{\bf 3. Linear Independence of Translates.}
In this section, we state the theorem on linear independence.
\bigskip

T{\eightrm HEOREM} 1.
{\it Let $\phi$ be a continuous, minimal (pre)scale function supported on
the interval $[0,N]$. Then the translates $\phi(\cdot - k),\,\, 
k = 0,\ldots,N-1$ are linearly independent over any set of positive
measure of the unit interval.}
\bigskip

{\bf Remarks.} As we noted above, this line of investigation was initiated
by Y.\ Meyer [10] and pursued by P.\ G.\ Lemari{\'e} in [7]. These
authors treated the case where the ``set of positive measure'' was the
entire unit interval.
Lemari{\'e} and Margouyres [8] gave another simplified proof that
showed the translates were linearly independent on any subinterval of
the unit interval.
Finally, Lemari{\'e} [7] showed that this property characterizes minimal
scale functions.
Those authors made no continuity assumptions.

{\bf Proof.}
We give a proof by contradiction as follows:  If the translates are
linearly dependent over a set of positive measure, we show that they are
dependent over a set of measure one in the unit interval. Since the function
$\phi$ is continuous, this means that the translates of $\phi$ are dependent
over the unit interval itself, thus contradicting the theorem of Meyer.

Throughout the proof, we will use  matrices $\bP_0$ and $\bP_1$.
To define these matrices, let us write the dilation equation for $\phi$
as
$$\phi(x/2) = \sum_{k=0}^N p_k \phi(x - k).$$
First, let us define the $(N-1)\times (N-1)$ matrix $\bP$ whose first row 
consists of the vector of odd numbered coefficients,
$p_{2k+1}$, followed by the 
appropriate number of zeros to give the vector $N-1$ components.  
The second row of $\bP$ is defined in the same way, using the even 
number coefficients, $p_{2k}$, followed by the appropriate number of zeros. 
Third and fourth rows are obtained from the first two rows by a cyclic
permutation of the indices: each entry is shifted to the right, with
the final entry, a zero, moving to first position. This procedure is
continued until $N-1$ rows are obtained. (Thus if $N = 2k$, the second
row will contain the $k+1$ entries $p_0,p_2,\ldots,p_{2k}$ followed by
$k-2$ zeros. The last row will contain $k-1$ zeros followed by the $k$
coefficients $p_1,p_3,\ldots,p_{2k-1}$.
If $N = 2k+1$, then the last row of the matrix consists of $k-1$ zeros,
followed by the $k+1$ entries $p_0,p_2,\ldots,p_{2k}$.)
Now define the two $N\times N$ matrices
$$\bP_0 = \pmatrix{p_0&\bp_t\cr
                     0&\bP\cr}
\qquad {\rm and}\qquad 
  \bP_1 = \pmatrix{\bP&0\cr
                     \bp_b&p_N\cr}$$
where $\bp_t$ is the $N-1$ vector consisting of the even numbered $p_k$,
starting with $p_2$, followed by the appropriate number of zeros;
$\bp_b$ consists of zeros followed by the coefficients $p_k$ where $k$
has the same parity as $N$, where the final entry of the vector $\bp_b$
is the coefficient $p_{N-2}$.

The roles of $\bP_0$ and $\bP_1$ are as follows: consider a general
linear combination of translates $\sum c_k \phi(x+k)$. If we take account
of the fact that $\phi$ is supported on $[0,N]$ and restrict attention
to $x\in [0,1]$, this sum is, in fact, finite and may be expressed as
$\sum_{k=0}^{N-1} c_k \phi(x+k)$. If we apply the dilation equation
to express each $\phi(\cdot + k)$ in terms of a sum of translates of
$\phi(2x)$, the resulting double sum is a certain linear combination
of translates of $\phi(2x)$ and $\phi(2x-1)$, depending on whether $x$
is in $[0, {1\over 2}]$ or in $[{1\over 2}, 1]$. The coefficients of this
linear combination are given by the matrices $\bP_0$ or $\bP_1$,
acting on the vector $\bc = (c_0,\ldots,c_{N-1})$.  (These matrices
are implicit in the reconstruction-decomposition schemes in the wavelet
literature, and appear explicitly, in the $3\times 3$ case in Daubechies
[4, section 7.2]. We summarize the above in the following proposition.
Let $\bPhi(x) = \big(\phi(x), \phi(x+1),\ldots,\phi(x + N-1)\big)^t$ for
$x\in [0,1]$, $\epsilon_k(x)$, $k = 1,2,\ldots$ be the digits in the
binary expansion of $x$. That is, $x = \sum\epsilon_i/2^k$, with
$\epsilon_k = 0$ or 1. Let $T$ be the plus-one shift on the
$\epsilon$-sequence: $T:(\epsilon_1,\epsilon_2,\ldots)\to
(\epsilon_2,\epsilon_3,\ldots)$.  We write $Tx = 
\sum_{k=1}^\infty\epsilon_{k+1}\big/ 2^k$.
\bigskip

P{\eightrm{ROPOSITION}}.
{\it For $\bc = (c_b,c_1,\ldots,c_{N-1})$ we have
$$\bc \circ \bPhi(x) = (\bP_{\epsilon_1} \bc^t)\circ \bPhi(Tx).$$
More generally, for any $m$,
$$\bc\circ\bPhi(x) 
      = (\bP_{\epsilon_1} \cdots \bP_{\epsilon_m}\bc^t)\circ\bPhi(T^m x).$$
}

\bigskip

{\bf Proof.}
Recall that the support of $\phi(\cdot + m)$ is the interval
$[-m, -m+N]$ in the following computation. For $x\in [0,1]$,
$$\eqalign{
\sum_{k=0}^{N-1} c_k\phi(x_k)
&= \sum_k c_k\bigg(\sum_j p_j\phi\big(2(x+k) - j\big)\bigg)\cr
&= \sum_m\bigg(\sum_k c_k p_{2k-m}\bigg)\phi(2x+m).\cr}$$
The inner sum is taken over all $k$ with the provision that $p_{2k-m} = 0$
if $2k-m$ is not one of the integers $0,1,\ldots,N-1$. The outer sum with
index $m$ changes according to whether $0 < 2x \le 1$ or $1 < 2x \le 2$,
due to the support condition mentioned above. In the first case, when
$\epsilon_1 = 0$, we have $0\le m\le N-1$; in the second case, when
$\epsilon_1 = 1$, $-1\le m\le N-2$. Thus, the transformation takes two
forms, with matrices $\bP_0$ and $\bP_1$. This proves the Proposition.

Now fix $\bc = (c_0,\ldots,c_{N-1})$ and let $K_\bc = \{x: \bc\circ\bPhi(x)
= 0\}$. The continuity of $\bPhi(x)$ implies that $K_\bc$ is closed. We
assume that $\bc\ne 0$, and that $K_\bc$ has positive measure in $[0,1]$;
we seek to contradict Meyer-Lemari{\'e} theorem. There are two cases to
consider.

{\it Case 1.}  There exists a finite sequence $\bP_{\epsilon_k},\,\,
k = 1,\ldots,m$ such that $\bP_{\epsilon_1}\cdots\bP_{\epsilon_m}\bc = 0$.
If this is the case, we have our contradiction since $\bc\circ\bPhi\equiv 0$
on the dyadic interval
$$\{x:\epsilon_1(x) = \epsilon_1,\ldots,\epsilon_m(x) = \epsilon_m\}.$$

{\it Case 2.}  The vector $\bc$ is such that $\bP_{\epsilon_1}\cdots
\bP_{\epsilon_m}\bc \ne 0$ for every finite sequence $\epsilon_1,\ldots,
\epsilon_m$.  In this case, we say that $\bc$ is a ``never zero'' vector.
\bigskip

L{\eightrm{EMMA}} 1. {\it Let $\bc$ be a never zero vector, and suppose
$m(K_\bc) > 0$. Then, for every $\eta$, $0 < \eta < 1$, there exists a
never zero vector $\bb$ such that $m(K_\bb) > 1 - \eta$.}
\bigskip

{\bf Proof.}
Since $K_\bc$ has positive measure, we can find dyadic interval
$I_j = \{x:  \epsilon_1(x) = \epsilon_1,\ldots,\epsilon_j(x) = \epsilon_j\}$
such that $m(K_\bc \cap I_j)/2^{-j} > 1 - \eta$.  This is a consequence of
the maximal martingale inequality, or alternatively, a point of density
argument.
Now apply the Proposition to points $x\in I_j$ to obtain a nonzero vector
$\bP_{\epsilon_1}\cdots\bP_{\epsilon_j}\bc = \bb$.  The set $K_\bb$ has
measure greater than $1-\eta$ since $m\big(\{x: x = T^j y\,\,\hbox{for some}
\,\, y \in I_j\}\big) = 1$.  This proves the Lemma.

Now set
$${\cal A}_\bc = \bigg\{ \ba \in \real^N:
\ba = {{\bb}\over {\Vert\bb\Vert_2}}\,\,{\rm for\,\, some\,\,}
\bb = \bP_{\epsilon_1}\cdots\bP_{\epsilon_j}\bc,\quad j\in \bZ\bigg\}.$$
By the Lemma, we have
$$\sup_{\ba\in {\cal A}_\bc} m(K_\bc) = 1.$$
Now we claim that the supremum is achieved: there is an
$\ba\in\bR^N$ with $\Vert\ba\Vert_2 = 1$ such that $m(K_\ba) = 1$.
To this end, we topologize the class ${\cal K}$ of compact sets
$K\subset [0,1]$ with a metric $\rho$ such that $({\cal K},\rho)$ is a
compact (Hausdorff) metric space. Let
$$\rho(A, B) = \sup_{x\in [0,1]} |d(x,A) - d(x,B)|$$
where $d(x,D) = \inf\{|x-y|: y\in D\}$. This metric is equivalent to
$$\sigma(A,B) = \inf\{\epsilon > 0, A\subset V_\epsilon(B)\,\,
{\rm and}\,\, B\subset V_\epsilon(A)\}$$
where $V_\epsilon(D) = \{z\in [0,1]: d(z,D) < \epsilon\}$.
Proofs of this equivalence and the fact that ${\cal K}$ is a compact space
may be found in Kornum [6, section 6.2].

Let $K_{\ba_n}$ be a sequence of sets such that $m(K_{\ba_n})$ tends to one.
>From this sequence, we may extract a convergent subsequence $K_{\ba_{n_k}}$.
Since $\Vert \ba_{n_k}\Vert_2 \equiv 1$, we may extract a convergent 
subsequence $\ba_{m_k}$, so that, finally, we obtain a sequence
$K_{\ba_n} \to K$ and $\ba_n \to \ba$. First of all, we claim that
$m(K) = 1$.  If not, there a $\delta > 0$ such that
$$m(K) \le m\big(V_j(K)\big) = d < 1.$$
In this case, there is an $n_0$ such that for all $n \ge n_0$,
$K_{\ba_n} \subset V_\delta(K)$, and so $m(K_{\ba_n}) \le m(V_\delta(K)) 
= d < 1$.
Since $m(K_{\ba_n}) \to 1$, this is a contradiction.
Second, we claim that $K\subset K_\ba$. Since $K_{\ba_n} \to K$, we have
$$\sup_{y\in [0,1]} |d(y,K) - d(y,K_{\ba_n})| \to 0.$$
Therefore, if $x\in K$, $d(x,K_{\ba_n}) \to 0$. 
Choose $x_n \in K_{\ba_n}$ so that
$|x - x_n| \to 0$. Then, by the Cauchy-Schwarz and triangle inequalities,
$$\eqalign{
|\ba \circ \Phi(x)|
&\le |(\ba - \ba_n) \circ \bPhi(x)|
     + |\ba_n \circ (\bPhi(x) - \bPhi(x_n))|
     + |\ba_n \circ \bPhi(x_n)|\cr
&\le \Vert \ba - \ba_n\Vert_2 \cdot \Vert\bPhi(x)\Vert_2
     + \Vert\bPhi(x) - \bPhi(x_n)\Vert_2.\cr}$$
Since both terms on the right tend to zero, we have the inclusion
$K\subset K_\ba$.
\bigskip

{\bf 4. Local Convergence of Wavelet Expansions.}
In [1], the following local convergence theorem is proved for Haar
series, using martingale methods.  

\bigskip

T{\eightrm{HEOREM}} A. {\it Let $f(x) = (f_0(x),f_1(x),\ldots)$
be a sequence of functions such that

{\rm (a)}  $f_j\in V_j$ where $\{V_k\}$ is the Haar multiresolution
analysis;

{\rm (b)}  $P_j(f_{j+1})(x) = f_j(x)$ for $j \ge 0$.

\noindent
Let $S^2(f)(x) = \sum (f_{j+1}(x) - f_j(x))^2 + f_0^2(x)$ and
$f^*(x) = \sup_j |f_j(x)|$.
Then, the following sets are equivalent almost everywhere:

{\rm (a)}  $\{x:  \lim_{j\to\infty} f_j(x)\,\, exists\,\, and\,\, is\,\,
finite\}$;

{\rm (b)}  $\{x:  S(f)(x) < + \infty\}$;

{\rm (c)}  $\{x:  \sup_j |f_j(x)| < \infty\}$.
}

\noindent
Gundy and Kazarian [2] extended this local convergence theorem to
the class of multiresolution analyses arising from the basic splines.
In fact, the proof did not appear to use properties specific to the
spline family. The basic regularity condition essential to the proof
is a two-norm condition, reminiscent of a condition first proposed by
Marcinkiewicz and Zygmund [9] in their study of series of
independent random variables. This condition, called condition (M-Z)
is as follows:

Let $\phi$ be a compactly supported scale function, supported on
$[0,N]$. We suppose that, for every $\delta$, $0 < \delta < 1$, there
exist constants $B_\delta$ and $C_\delta$ such that

(M-Z)  For every measurable subset $E\subset [0,1]$ of measure greater
than $\delta$, and any sequence $a_k$; $k = 0,1,\ldots,N-1$, we have
$$\eqalign{
C_\delta \sum_{k=0}^{N-1} |a_k|
&\le \sup_{x\in E} \bigg|\sum_{k=0}^{N-1} a_k \phi(x+k)\bigg|\cr
&\le B_\delta \sum_{k=0}^{N-1} |a_k|.\cr}$$
The constants $B_\delta, C_\delta$ depend only on $\phi$ and the measure
of the set $E$.  The condition holds for the class of $B$-spline scale
functions, as pointed out in [2]. However, the scope of the condition
was not known, and left as an open problem in [2]. The following
theorem answers this question.

\bigskip

L{\eightrm{EMMA}} 2.  {\it Let $\{V_j\}$ be a multiresolution analysis
such that $V_0$ contains continuous functions of compact support.  Then
the minimal scale function $\phi$ satisfies condition (M-Z).}
\bigskip

Before proving the lemma, we state the following theorem, in which we use
the definitions in Theorem A, for multiresolution analyses more general
than the Haar system.
\bigskip

T{\eightrm{HEOREM}} 2. (Theorem B of [2])
{\it Let $\{V_j\}$ be a multiresolution analysis that contains
continuous functions of compact support. Then the following sets are
equivalent almost everywhere:

{\rm (a)}  $\{x:  \lim_{j\to\infty} f_j(x)\,\, exists\,\, and\,\, is\,\,
finite\}$;

{\rm (b)}  $\{x:  S(f)(x) < + \infty\}$;

{\rm (c)}  $\{x:  \sup_j |f_j(x)| < \infty\}$.
}

\bigskip

{\bf Proof of Lemma 2.}
Since $\phi$ is continuous on $[0,1]$, the issue is to show the existence
of $C_\delta$ that is uniform over all sets $E\subset [0,1]$ of measure
greater than $\delta$.  First, observe that, since the translates of $\phi$
are linearly independent over $E$ (Theorem 1), there is a constant $C(E)$
such that
$$C(E) \sum_{k=0}^{N-1} |a_k|
\le \sup_{x\in E} \bigg|\sum_{k=0}^{N-1} a_k \phi(x+k)\bigg|.$$
This follows from the fact that the right-hand side defines a norm on
$\bR^N$: the linear independence of the translates of $\phi$ guarantees
that this norm is strictly positive on $\bR^N - \{0\}$. Since the
left-hand side is also a norm, the existence of a constant is assured by the
open mapping theorem.
Now we must show that
$$\inf \{C(E): m(E) \ge \delta\} > 0.$$
It is enough to show this for closed sets. To this end, we show that
$$C(E) = \sup_{x\in E} {{|\ba\circ\bPhi(x)|}\over {\Vert\ba\Vert_2}}$$
is a continuous function of $E$ for the topology $({\cal K}, \rho)$
introduced above. Let $\epsilon > 0$ be given, and let $\{A_n\}$ be a
sequence of sets converging to $A$ in ${\cal K}$. The vector $\bPhi(x)$
is uniformly continuous on $[0,1]$, so that
$$\Vert\bPhi(x) - \bPhi(y)\Vert_2
\le {{\epsilon\Vert\ba\Vert_1}\over {\Vert\ba\Vert_2}}$$
whenever $|x-y| \le \delta(\epsilon)$. Let $n_0$ be an integer such that
$$A_n\subset V_\delta(A)\qquad {\rm and}\qquad
A_n\subset V_\delta(A_n)$$
for all $n \ge n_0$. Since $\bPhi$ is continuous on $\real^N$ and
$\overline{V_\delta(A_n)}$ is compact, it follows that for $n\ge n_0$,
there exists an $x_n\in \overline{V_\delta(A_n)}$ for which
$$C(A) \le C(\overline{V_\delta(A_n)})
= {{|\ba\circ\bPhi(x_n)|}\over {\Vert\ba_1\Vert}}.$$
By the uniform continuity of $\bPhi(\cdot)$, 
$|C(A_n) - C(\overline{V_\delta(A_n)})| \le \epsilon$, so that
$$C(A) \le C(A_n) + \epsilon.$$
If we reverse the roles of $A_n$ and $A$ in the above argument, we see that
$C(A_n) \le C(A) + \epsilon$. Thus, $C(E)$ is continuous on ${\cal K}$.
The collection $\{E\in {\cal K}: m(E) \ge\delta\}$ is closed in ${\cal K}$
by a similar argument, so there is a closed $E_0$, $m(E_0)\ge \delta$
such that $0 < C(E_0) \le C(E)$ for any $E$, $m(E) \ge\delta$.

\bigskip

{\bf 5. Concluding Remarks.}
The quadratic variation functional $S(f)$ of Theorem 3 is invariant
under changes of scale functions $\phi$ for $V_0$ and $\psi$ for $V_1$:
$S(f)$ is defined from the sequence of projections $\{P_j\}$
without specific reference to the choice of scale function. However,
$S(f)$ is an ``incomplete'' square function in the sense that if the
prewavelet family $\{\psi(2^j x - k)\}$ is orthogonalized in $k$,
to obtain a family $\{\tilde\psi (2^j x-k)\}$ that is orthonormal
in both variables $j,k\in \bZ$, then one could consider a quadratic
variation functional $\bS^2(f)(x) 
= \sum_{j,k} (a_{j,k}\tilde\psi (2^j x-k))^2$.  If we have a multiresolution
analysis that admits a compactly supported, continuous orthonormal family
$\{\tilde\psi (2^j x-k)\}$, then one can show that $S(f)(x)$
and $\bS(f)(x)$ are finite on the same set, up to a set of measure
zero. The proof of this fact follows the same lines as the proof in
[2].  Since the details are given there, we will not repeat them here.
\vfil\eject

\centerline{Bibliography}
\bigskip
\bigskip

1. R.\ F.\ Gundy, Martingale theory and the pointwise convergence of certain
orthogonal series, {\it Trans.\ Amer.\ Math.\ Soc.} {\bf 124} (2) (1966),
228-248.

2. R.\ F.\ Gundy, K.\ Kazarian, Stopping times and local convergence for
spline wavelet expansions, Rutgers University Technical Report \#97-002.

3. I.\ Daubechies, Orthonormal bases of compactly supported wavelets,
{\it Comm.\ Pure Appl.\ Math.} {\bf 41} (1988), 909-996.

4. I.\ Daubechies, {\it Ten Lectures on Wavelets}, CBS-NSF Regional
Conference in Applied Math.,  {\it SIAM} {\bf 61} (1992).

5. E.\ Hern{\'a}ndez, G.\ Weiss, {\it A First Course on Wavelets},
CRC Press, Boca Raton FL, 1996.

6. P.\ Kornum, {\it Construction of Borel Measures on Metric Spaces},
Aarhus Universitet, Matematisk Institut, 1980.

7. P.\ G.\ Lemari{\'e}, Fonctions {\`a} support compact dans les analyses
multi-r{\'e}solutions, {\it Rev.\ Mat.\ Iberoamericana} {\bf 7} (2) (1991),
157-182.

8. P.\ G.\ Lemari{\'e}, G.\ Malgouyres, Support des fonctions de base dans
une analyse multi-r{\'e}solution, {\it C.\ R.\ Acad.\ Sci.\ Paris},
S{\'e}ries I, {\bf 313}, 377-380.

9.  J.\ Marcinkiewicz, A.\ Zygmund, Sur les fonctions ind{\'e}pendentes,
{\it Fund.\ Math.} {\bf 29} (1937), 60-90.

10. Y.\ Meyer, Ondelettes sur l'intervalle, {\it R.\ Mat.\
Iberoamericana} {\bf 7} (2) (1991), 115-133.

\end